\documentclass[11pt]{amsart}
\usepackage{amssymb, latexsym}

\newcommand{\ch}{\mbox {\bf 1}}
\newcommand{\ex}{{\mathbb {E}}}

\newcommand{\witi}{\widetilde}
\newcommand{\R}{{\mathbb{R}}} %

\newcommand{\Z}{{\mathbb{Z}}}
\newcommand{\N}{{\mathbb{N}}}

\newcommand{\ol}[1]{\overline{#1}}

\newcommand{\re}{\mbox{\rm Re}}

\newtheorem{itdefinition}{Definition}

\newtheorem*{itdefinition*}{Definition}
\newenvironment{definition*}{\begin{itdefinition*}\rm}{\end{itdefinition*}}
\newtheorem{assumption}{Assumption}

\newtheorem{theorem}{Theorem}
\newtheorem*{theorem*}{Theorem}
\newtheorem{itremark}{Remark}

\newtheorem*{itremark*}{Remark}
\newtheorem*{corollary*}{Corollary}
\newenvironment{remark*}{\begin{itremark*}\rm}{\end{itremark*}}
\newtheorem{prop}{Proposition}

\newtheorem{itexample}{Example}

\begin{document}
\title[ Elliptic Operators with Nonlocal Boundary
Condition] { Spectral Analysis of a Family of Second-Order
Elliptic Operators with Nonlocal Boundary Condition Indexed by a
Probabilty Measure}

\author{Iddo Ben-Ari}
\address{Department of Mathematics\\
University of California - Irvine\\
Irvine, CA 92697\\
        USA}
       \email{ibenari@math.uci.edu}

\author{Ross G. Pinsky}
\address{Department of Mathematics\\
Technion---Israel Institute of Technology\\
Haifa, 32000\\
Israel} \email{pinsky@math.technion.ac.il} \subjclass[2000]{}
\date{}

 \begin{abstract} Let $D\subset R^d$ be a bounded domain  and let
\[
L=\frac12\nabla\cdot a\nabla +b\cdot\nabla
\]
be a second order elliptic operator on $D$.
Let $\nu$ be a probability measure on $D$.
   Denote by ${\mathcal L}$ the differential operator
  whose domain is specified by the following  non-local boundary condition:
 $$
 {\mathcal D_{{\mathcal   L}}}=\{f\in C^2(\ol{D}): \int_D f d\nu = f|_{\partial D}\},
 $$
 and which coincides with $L$ on its domain.
Clearly 0 is an eigenvalue for
 $\mathcal L$, with the  corresponding
eigenfunction being constant. It is known that
$\mathcal L$ possesses an infinite sequence of eigenvalues, and that with the exception
of the zero eigenvalue, all eigenvalues have negative real part.
Define the spectral gap of $\mathcal {L}$, indexed by $\nu$,  by
\begin{equation*}
\gamma_1(\nu)\equiv\sup\{\re~\lambda:0\neq \lambda\mbox{ is an
eigenvalue for }{\mathcal L}\}.
\end{equation*}
In this paper we investigate the eigenvalues of $\mathcal L$ in
general and the spectral gap $\gamma_1(\nu)$ in particular.

The
operator $\mathcal L$ and its spectral gap $\gamma_1(\nu)$ have
probabilistic significance.
The operator $\mathcal L$ is the generator of a diffusion process with random jumps from the
boundary, and $\gamma_1(\nu)$ measures the exponential rate of convergence of this process
to its invariant  measure.
\end{abstract}
\maketitle

\setcounter{equation}{0}\section{Introduction and Statement of
Results} Let $D\subset R^d$ be a bounded domain  and let
\[
L=\frac12\nabla\cdot a\nabla +b\cdot\nabla
\]
be a second order elliptic operator on $D$. We will assume that
 $a=\{a_{ij}\}_{i,j=1}^{d}$ is positive definite with entries in $C^{2,\alpha}(\R^d)$  and that
$b=(b_1,\dots,b_d)$ has entries in $C^{1,\alpha}(\R^d)$,
for some $\alpha\in(0,1]$.  Note that we have written the
principal part of the operator $L$ in divergence form. This has
been done  for convenience and, in light of the above conditions
on the coefficients,  without loss of generality.
We will assume either that $D$ has a $C^{2,\alpha}$-boundary or that $D=D_1\times\cdots\times D_k$, and
$L=\sum_{i=1}^k L_i$, where $L_i$ is defined on $D_i$ and
$D_i$ has a $C^{2,\alpha}$-boundary. This latter situation allows in particular for the case of  $\frac12\Delta$ on a cube.

Let $\nu$ be a probability measure on $D$.
   Denote by ${\mathcal L}$ the differential operator
  whose domain is specified by a non-local boundary condition as follows:
 $$
 {\mathcal D_{{\mathcal   L}}}=\{f\in C^2(\ol{D}): \int_D f d\nu = f|_{\partial D}\},
 $$
 and which coincides with $L$ on its domain.
(Non-local boundary conditions in the spirit of the one above in the context of    parabolic operators can be found in the physics
literature on ``well-stirred'' liquids. See \cite{vandenberg} and \cite{CJ}.)

Clearly 0 is an eigenvalue for
 $\mathcal L$, with the  corresponding
eigenfunction being constant. It is known that
$\mathcal L$ possesses an infinite sequence of eigenvalues, and that with the exception
of the zero eigenvalue, all eigenvalues have negative real part (see Theorem BP below).
Note that the
  operator $\mathcal L$ depends on the measure $\nu$ through its domain of definition.  Define
 the spectral gap of $\mathcal {L}$, indexed by $\nu$,  by
\begin{equation}\label{gap}
\gamma_1(\nu)\equiv\sup\{\re~\lambda:0\neq \lambda\mbox{ is an
eigenvalue for }{\mathcal L}\}.
\end{equation}
In this paper we investigate the eigenvalues of $\mathcal L$ in
general and the spectral gap $\gamma_1(\nu)$ in particular. The
operator $\mathcal L$ and its spectral gap $\gamma_1(\nu)$ have
probabilistic significance which we now point out.

Let $G^D(x,y)$ denote the Green's function for $L$,   defined by
  $$G^D(x,y)=\int_0^{\infty}p^D(t,x,y)dt,$$
where $p^D(t,x,y)$ is the Dirichlet heat kernel for $L-\frac
d{\partial t}$ in $D$, or equivalently, as a function of $y$,
$p^D(t,x,y)$ is the transition subprobability density for the
diffusion process $Y(t)$ in $D$ corresponding to $L$, starting
from $x\in D$ and killed upon exiting $D$. It was shown in
\cite{BP} that there exists a Markov  process $X(t)$ in $D$ which
coincides with the diffusion $Y(t)$ governed by $L$ until it exits
$D$, at which time it jumps to  a point in the domain according to
the distribution $\nu$ and starts the diffusion afresh. This same
mechanism is repeated independently each time the process reaches
the boundary. This process is   called a \it diffusion with
random jumps from the boundary.\rm\ In light of the above
probabilistic connection, from now on we will refer to the measure
$\nu$ appearing in the definition of $\mathcal L$ as the jump
measure. Denote expected values corresponding to this process
starting from $x\in D$ by $\ex_x$. Let ${\mathcal P}(D)$ denote
the space of probability measures on $D$.
Under the smoothness conditions stated above,
 the
following theorem was proven in \cite[Theorem 1 and the remark following it]{BP}.
\begin{theorem*}[BP]
 \label{ergodic}
 Let $X$ be the diffusion with random jumps from the boundary
corresponding to $L$ and $\nu$.\\
\noindent i. There exists a unique invariant measure $\mu$ for the process. It has a density,  also denoted by $\mu$, which is given by
$$
\mu(y)=\frac{\int_DG^D(x,y)d\nu(x)}{\int_D\int_DG^D(x,z)d\nu(x)dz}.
$$
The map  Inv:\ ${\mathcal P}(D) \rightarrow {\mathcal P}(D)$, defined by Inv$(\nu)=\mu$, is continuous in the topology of weak convergence of probability measures.\\
\noindent ii. The operator $\mathcal L$ possesses an infinite sequence of eigenvalues.
Furthermore,
$$\lim_{t\to\infty}\frac1t\log\sup_{f\in L^{\infty}(D),~||f||_\infty\le1}||\ex_x f(X(t))-
\int_D fd\mu||_\infty=\gamma_1(\nu)<0,
$$
where $\gamma_1(\nu)$, defined in \eqref{gap}, is the spectral gap of $\mathcal L$.
\end{theorem*}

\begin{remark*}
Actually, part (ii) of  Theorem BP was proved in \cite{BP} for  a more general problem, where the jump measure
from the boundary is allowed to depend on the boundary location.
\end{remark*}

We now turn to the analysis of the eigenvalues of $\mathcal L$ in
general and of the spectral gap of $\mathcal L$ in particular.
Note that by Theorem BP, the larger the spectral gap, the faster
is the rate of convergence to equilibrium for the diffusion with
random jumps.

We begin with a very special case of jump measure $\nu$ where the eigenvalues (and eigenfunctions)
of $\mathcal L$ can be completely characterized in terms of those of $L$ with the Dirichlet
boundary condition.
Recall that the operator $L$ with the Dirichlet boundary condition possesses
an infinite sequence of eigenvalues, all of which have negative real part.
By the Krein-Rutman theorem, the \it principal eigenvalue\rm\---the eigenvalue with largest real part---is
real and simple, and the corresponding eigenfunction does not change sign
The same is true for  $\tilde L$, the formal adjoint of $L$ with the Dirichlet
boundary condition. Furthermore, the principal eigenvalues of $L$ and $\tilde L$ coincide.
Let $\witi{\phi^D_0}>0$ denote the principal eigenfunction
corresponding to the principal eigenvalue for  $\witi{L}$.
 Normalize it by $\int_D\witi{\phi^D_0}(x)dx=1$. Abusing
notation, we will also let $\witi{\phi^D_0}$ denote the measure
with density $\witi{\phi^D_0}$. (The measure $\witi{\phi^D_0}$ is
the so-called quasi-invariant distribution for the original
diffusion corresponding to $L$ with killing at the boundary. That
is, one has $\ex^D_{\witi{\phi^D_0}}(f(Y(t))|\tau_D>t)=\int_D f(x)
\witi{ \phi^D_0}(x) dx$, for all $t>0$, where  $\tau_D$ is the
first exit time of the diffusion $Y(t)$ from $D$ and
$\ex^D_{\witi{\phi^D_0}}$ denotes the expectation for the
diffusion killed at the boundary and starting from the
distribution  $\witi{\phi^D_0}$.)
\begin{theorem}~
\label{thm:quasi inv}
Consider the operator $\mathcal L$ in the case that the jump
measure is given by $\nu=\witi{\phi^D_0}$, where $\witi{\phi^D_0}$ is
the normalized principal eigenfunction for the formal adjoint $\tilde L$ of $L$ with the Dirichlet boundary
condition.
Let $\{\lambda_n^D\}_{n=0}^\infty$ denote the eigenvalues for $L$
with the Dirichlet boundary condition, labeled so that Re\
$\lambda^D_{n+1}\le$ Re\ $\lambda^D_n$, and let
$\{\phi_n\}_{n=0}^\infty$ be a corresponding sequence of
eigenfunctions. Then the eigenvalues for ${\mathcal L}$ are 0 and
$\{\lambda^D_n\}_{n=1}^\infty$ and a  corresponding sequence of
eigenfunctions is given by 1 and $\{\phi_n\}_{n=1}^\infty$. In
particular,
$$
\gamma_1(\witi{\phi^D_0})=\text{Re} \ (\lambda^D_1).
$$
Furthermore, $\witi{\phi^D_0}$ is the invariant probability
measure for the diffusion with random jumps from the boundary
corresponding to $\mathcal L$. In fact, $\witi{\phi^D_0}$ is the
unique fixed point for the map Inv:\ ${\mathcal
P}(D)\rightarrow{\mathcal P}(D)$ defined in Theorem BP.
\end{theorem}


In order to make the spectral analysis tractable
when the jump measure $\nu$ is not the special measure considered in Theorem \ref{thm:quasi inv},
we will need to assume that
the operator $L$ with the Dirichlet boundary condition is self-adjoint,
although $\mathcal L$ will still not be self-adjoint, as we now explain.
If the first-order term  $b$ in  the operator $L$
is of the form $b=a\nabla Q$, then the operator $L$
can be written in the form $L=\frac12\exp(-2Q)\nabla\cdot a\exp(2Q)\nabla$.
Since we can replace $Q$ by $Q+c$, where $c$
is a constant, without changing $L$,
we will assume without loss of generality  that $\int_D\exp(2Q)dx=1$.
Let $\mu_{\text{rev}}$ denote the probability measure $\exp(2Q)dx$.
In this case, the operator  $L$
 with the Dirichlet boundary condition, considered as an operator on
 $L^2(D,\mu_{\text{rev}})$,
is symmetric on the domain of smooth functions vanishing at the boundary and is
 self-adjoint on an appropriate domain of definition.
(The diffusion process in $D$ killed at the boundary,
corresponding to $L$ with the Dirichlet boundary condition, is reversible and the normalized
reversible
measure is $\mu_{\text{rev}}$; whence the notation $\mu_{\text{rev}}$.)
The operator ${\mathcal L}$, on the other hand, will \it never\rm\
be self-adjoint. Indeed,  a straight forward calculation shows
that the adjoint operator (with respect to Lebesgue measure)
$\tilde{\mathcal L}$ of $\mathcal L$ is defined on a domain which
includes $\{v\in C^2(D)\cap C(\bar D):v=0 \ \text{on}\ \partial
D\}$, and for such functions one has $\tilde {\mathcal L}v=\tilde
Lv-(\int_D\tilde Lv)\nu$, where $\tilde L=\frac12\nabla \cdot
a\nabla-b\nabla-\nabla\cdot b$. In the case that $L$ is
self-adjoint, if one takes the adjoint  of $\mathcal L$ with
respect to $\mu_{\text{rev}}$, then the adjoint is defined on the
above class of functions by $\tilde{\mathcal L}v=Lv-(\int_D
Lv)\nu$.

We will begin with a key theoretical result, which will be mined to obtain more concrete results.
Before we can state the theorem, we need some additional notation.
The eigenvalues of the self adjoint operator $L$ are real and
negative. We will denote them by $\{\lambda_n^D\}_{n=0}^\infty$,
labelled in nonincreasing order.
 Denote the corresponding eigenfunctions by $\{\phi^D_n\}_{n=0}^\infty$, normalized by
$\int_D\phi^D_nd\mu_{\text{rev}}=1$, $n\ge0$, and $\phi_0^D>0$.
Let
\begin{equation}\label{FG}
F_n\equiv\int_D\phi^D_nd\mu_{\text{rev}} \ \ \text{and}\ \ G_n(\nu)\equiv\int_D\phi_n^Dd\nu.
\end{equation}
Let $\{\Lambda^D_n\}_{n=0}^\infty$ denote the collection of distinct eigenvalues among $\{\lambda^D_n\}_{n=0}^\infty$, labelled in decreasing order.
We will sometimes need the following assumption.
\begin{assumption}
The Fourier series
$\sum_{n=0}^\infty\frac{F_n}{\lambda_n^D}\phi_n^D(x)$ converges
uniformly and absolutely.
\end{assumption}

\begin{theorem}\label{E}
Assume that
the operator  $L$ with the Dirichlet boundary condition is self-adjoint.
Let
\[
E_\nu(\lambda)\equiv\sum_{n=0}^\infty\frac{F_nG_n(\nu)}{\lambda_n^D-\lambda}.
\]
Let $d_n$ denote
the dimension of the eigenspace corresponding to the $n$-th distinct eigenvalue $\Lambda_n^D$ of $L$.
Assume either that $\nu$ possesses an $L^2(D,d\mu_{\text{rev}})$-density
or that Assumption 1 on the operator $L$ holds.

Then the  set of nonzero eigenvalues of $\mathcal{L}$  and their multiplicities are given as follows:

{\bf\noindent i.}
The set $\{\lambda:E_\nu(\lambda)=0\}-  \{\Lambda^D_n\}_{n=1}^\infty$ consists of simple eigenvalues;

{\bf\noindent ii.} For each $n=1,2\cdots$, the following rule determines whether $\Lambda_n^D$ is an eigenvalue, and if so,
specifies its multiplicity:

\noindent
If $d_n=1$ and  neither $F_m=0$ nor $G_m(\nu)=0$, for the $m$
satisfying $\lambda_m^D=\Lambda_n^D$, then $\Lambda_n^D$ is not an eigenvalue.
 Otherwise, $\Lambda_n^D$ is an eigenvalue and its multiplicity is specified as follows:

\noindent If $G_m(\nu)\neq0$ for some $m$ such that $\lambda_m^D=\Lambda^D_n$
and $F_m\neq0$ for some $m$ such that $\lambda_m^D=\Lambda^D_n$, then
the multiplicity is $d_n-1$;

\noindent If $G_m(\nu)=0$ for all $m$ such that $\lambda_m^D=\Lambda^D_n$
and $F_m\neq0$ for some $m$ such that $\lambda_m^D=\Lambda^D_n$, or if
 $G_m(\nu)\neq0$ for some $m$ such that
$\lambda_m^D=\Lambda^D_n$ and
$F_m=0$ for all $m$ such that $\lambda_m^D=\Lambda^D_n$, then the multiplicity
is $d_n$;

\noindent
If $G_m(\nu)=0$ for all $m$ such that
$\lambda_m^D=\Lambda^D_n$ and
$F_m=0$ for all $m$ such that $\lambda_m^D=\Lambda^D_n$,
then the multiplicity is $d_n$ if
$E_\nu(\Lambda_n^D)\neq0$ and is $d_{n+1}$ if
$E_\nu(\Lambda_n^D)=0$.

Furthermore, even without Assumption 1 or the density condition on $\nu$,  the set of eigenvalues of $\mathcal{L}$ includes
those listed in (ii).

\end{theorem}

Note that the complete characterization of the spectrum in
Theorem \ref{E} always holds
if the jump measure $\nu$ possesses an $L^2(D,\mu_{\text{rev}})$-density.
If the operator $L$ on $D$ satisfies Assumption 1, then it holds for all jump measures $\nu\in\mathcal P(D)$.
The following theorem collects some  sufficient conditions for Assumption 1 to hold.

\begin{theorem}\label{conditions}
i. If $d=1$, then Assumption 1  holds.

\noindent ii. Let
$d=2$ and let $L=\frac12\exp(-2Q)\nabla\cdot a\exp(2Q)\nabla$
satisfy $Q=\frac12\log\sqrt{det(a^{-1}})$ (in which case $L$ can
be considered as $\frac12\Delta_M$, where $\Delta_M$ is the
Laplacian of a Riemannian manifold with metric $a$). Then
Assumption 1  holds.

\noindent iii If $d\le 3$ and the eigenfunctions $\{\phi_n^D\}_{n=0}^\infty$
are uniformly bounded, then Assumption 1  holds.
\end{theorem}

\begin{remark*}
A direct calculation (see the proof of Proposition \ref{cube} below) shows that
the eigenfunctions $\{\phi_n^D\}_{n=0}^\infty$ are uniformly bounded for $L=\frac12\Delta$
in $D=(0,1)^d$; however such a bound does not  hold if $D$ is a sphere \cite{Gri02}.
\end{remark*}

As a first application of Theorem \ref{E}, we identify a class of
jump measures $\nu$ for which all the eigenvalues of $\mathcal L$
are real. The analysis of the spectrum in this case turns out to
be more tractable.

\begin{theorem}\label{real}
Assume that
the conditions of Theorem \ref{E} are in force,
and let $F_n$ and $G_n(\nu)$ be as in \eqref{FG}.
Assume also that the jump measure $\nu$ satisifies
one of the following two conditions:

i.  $F_nG_n(\nu)\ge0$, for all $n\ge1$, or $F_nG_n(\nu)\le0$, for all $n\ge1$;

ii. $F_nG_n\neq0$ for at most two values of $n$.

\noindent Then all the eigenvalues of  $\mathcal L$ are real.
\end{theorem}

When the nonzero eigenvalue with the largest real part is real,  we can prove an upper bound on
the eigenvalue spectral gap, $\gamma_1(\nu)$, of $\mathcal L$.

\begin{theorem}\label{upperbound}
If the jump measure  $\nu$ is such that the nonzero eigenvalue of
$\mathcal L$ with the largest real part is real, then
$$
\gamma_1(\nu)<\lambda_0^D.
$$
\end{theorem}

\bf\noindent Remark 1.\rm\
Theorem \ref{upperbound} holds regardless of whether the operator $L$ is self-adjoint; however,
if $L$ is not self-adjoint then we have no way of determining whether the nonzero eigenvalue
of $\mathcal L$ with the largest real part is real.

\noindent \bf Remark 2.\rm\  As is well known,
$\lambda_0^D$ gives the exponential rate of decay in $t$ of the probability
that the diffusion $Y(t)$ in $D$ corresponding to $L$
has not yet
hit the boundary by time $t$; that is, $\lim_{t\to\infty}\frac1t\log P_x(\tau_D>t)=\lambda_0^D$,
where $\tau_D$ is the first exit time of the diffusion from $D$.
Now since
 $\gamma_1(\nu)$ gives the exponential rate of convergence
of the distribution of the diffusion with random jumps to its
invariant measure, and since the jump mechanism only comes into affect after time
$\tau_D$, Theorem \ref{upperbound} might seem (at least at first blush) counter-intuitive.


The
 normalized reversible measure $\mu_{rev}$, with respect to which $L$ is self-adjoint,
   plays a distinguished role as the jump measure.
In particular, in this case the spectral gap can be given by a variational formula.

\begin{theorem}\label{Ereversible}
Assume that
the operator  $L$ with the Dirichlet boundary condition is self-adjoint.
Let $\{\lambda_n^D\}_{n=0}^\infty$ denote the eigenvalues of $L$ with the Dirichlet boundary condition.
Let the jump measure
be the normalized reversible measure $\mu_{rev}$.
\begin{enumerate}
\item[i.] All the eigenvalues of $\mathcal L$ are real.
\item[ii.]
$$
\gamma_1(\mu_{\text{rev}})=-\inf\frac{\frac12\int_D(\nabla ua\nabla u)d\mu_{\text{rev}}}{\int_Du^2d\mu_{\text{rev}}},
$$
where the infimum is over functions $u\neq0$ satisfying
$u|_{\partial D}=\int_Dud\mu_{\text{rev}}=0$. The infimum is
attained at a function $u_{min}$ which satisfies the equation
$Lu=\gamma_1(\mu_{\text{rev}})u+C$, for some constant $C$, and the
eigenfunction $v_1$ for $\mathcal{L}$ corresponding to the
eigenvalue $\gamma_1(\mu_{\text{rev}})$ is given by
$v_1=u_{min}+\frac C{\gamma_1(\mu_{\text{rev}})}$.
\item[iii.]
\[
\lambda_1^D\le\gamma_1(\mu_{rev})<\lambda_0^D.
\]
More precisely, consider the function
\[
E_{\mu_{rev}}(\lambda)=\sum_{n=0}^\infty\frac{F^2_n}{\lambda_n^D-\lambda},
\]
which is increasing for $\lambda\in(\lambda_1^D,\lambda_0^D)$.
If the equation $E_{\mu_{rev}}(\lambda)=0$ possesses a root in $(\lambda_1^D,\lambda_0^D)$, then
$\gamma_1(\mu_{\text{rev}})$ is equal to this root. Otherwise,
$\gamma_1(\mu_{\text{rev}})=\lambda_1^D$.
In particular, such a root will exist
if $F_j=\int_D\phi^D_jd\mu_{\text{rev}}\neq0$, for some $j\in\{1,\cdots,k_0\}$,
where $k_0=max\{n:\lambda_n^D=\lambda_1^D\}$.
 If $F_j=\int_D\phi^D_jd\mu_{\text{rev}}=0$, for all $j\in\{1,\cdots,k_0\}$,  then
$\gamma_1(\mu_{\text{rev}})>\lambda_1^D$ if and only if
\begin{equation}\label{borderline}
\frac{F_0^2}{\lambda_0^D-\lambda_1^D}<\sum_{n=k_0+1}^\infty\frac{F^2_n}{\lambda_1^D-\lambda_n^D}.
\end{equation}
\end{enumerate}

\end{theorem}

\begin{remark*}
Consider the diffusion process corresponding to $L$ as in Theorem \ref{Ereversible}
with \it reflection \rm\ at the boundary in the conormal direction $an$, where
$n$ denotes the inward unit normal to $D$. The process is reversible and
it  corresponds to a self-adjoint operator on $L^2(D,\mu_{\text{rev}})$ which is
an extension of  $L$ with the Neumann boundary condition
$\nabla u \cdot an=0$ on $\partial D$.
For this process, $\mu_{\text{rev}}$ is the invariant measure, and the rate of convergence to
$\mu_{\text{rev}}$  is given by the largest nonzero eigenvalue, $\lambda_1^N$.
This eigenvalue is given by the variational formula in part (ii)
of Theorem \ref{Ereversible}, but with the infimum being taken
over functions $u$ satisfying $\int_Dud\mu_{\text{rev}}=0$ (without the additional restriction
that $u|_{\partial D}=0$).
The infimum is attained  at the eigenfunction(s) corresponding to $\lambda_1^N$, and it is known
that any such function does not vanish identically on $\partial D$. Thus, it follows from part (ii) of Theorem
\ref{Ereversible} that $\lambda_1^N>\gamma_1(\mu_{\text{rev}})$. Therefore, the rate of convergence to
equilibrium is greater for the diffusion with random jumps with jump measure
 $\mu_{\text{rev}}$ than for the reflected diffusion, whose invariant measure is $\mu_{\text{rev}}$.
\end{remark*}

Here is an application of  condition \eqref{borderline} in part (iii) of Theorem
\ref{Ereversible}.

\begin{prop}\label{cube}
Consider the operator $\frac12\Delta$ in
the $d$-dimensional unit cube, $D=(0,1)^d$, and let the jump measure be Lebesgue
measure, $l_d$, on $D$. One has $\lambda_0^D=-\frac{d\pi^2}2$ and
$\lambda_1^D=-\frac{(d+3)\pi^2}2$.
\begin{enumerate}
\item [i]
If $d\le10$, then $\gamma_1(l_d)=\lambda_1^D$.
\item [ii.] If $d\ge11$, then
$\lambda_1^D<\gamma_1(l_d)<\lambda_0^D$.
\end{enumerate}
\end{prop}
\begin{remark*}
Note that  $\gamma_1(l_d)$ decreases to $-\infty$
as $d\to\infty$. Thus, for Brownian motion
 in the $d$-dimensional
cube with random jumps from the boundary
with normalized Lebesgue measure as the jump measure,
the rate of convergence to equilibrium becomes arbitrarily
fast as the dimension increases. This is because starting from any point, the distribution
of the hitting time of the boundary converges to the $\delta$-measure at 0
as $d\to\infty$, which means that as $d\to\infty$, the process constantly gets redistributed according to Lebesgue
measure after  arbitrarily small intervals of time.
In contrast to this, consider  Brownian motion in the $d$-dimensional unit cube
with normal reflection at the boundary. The rate of convergence
to equilibrium is governed by the largest  nonzero eigenvalue
of the Neumann Laplacian, which  is $\lambda_1^N=-\frac{\pi^2}2$, independent
of $d$. Similarly, consider  Brownian motion in the $d$-dimensional unit cube,
conditioned never to hit the boundary \cite{Pin85}. This process corresponds to the $h$-transformed
operator $(\frac12\Delta-\lambda^D_0)^{\phi^D_0}$.
The rate of convergence to equilibrium
is governed by the largest nonzero eigenvalue, which is
 $\lambda_1^D-\lambda_0^D=-\frac{(d+3)\pi^2}2+\frac{d\pi^2}2=-\frac32\pi^2$, independent of $d$.
\end{remark*}

We have the following result for the one-dimension Laplacian.

\begin{prop}\label{1dimBM}
Consider the operator $\frac12\frac{d^2}{dx^2}$ in the interval $(0,1)$.

\begin{enumerate}
\item[i.] If the jump measure is deterministic; that is,
 $\nu=\delta_p$, for some $p\in(0,1)$, then $\gamma_1(\nu)=\lambda_1^D=-2\pi^2$;

\item[ii.]
If the jump measure $\nu$ is such that the nonzero eigenvalue of
$\mathcal L$ with the largest real part is real, then
$$
\gamma_1(\nu)=\lambda_1^D=-2\pi^2.
$$
\end{enumerate}
\end{prop}

\begin{remark*}
Part (i) above was shown in \cite{BM_8} and \cite{BM_reb}.
(Actually, $-\frac{\pi^2}2$ was obtained in \cite{BM_8}, because a
certain cancellation was not taken into account. The correct
result appears in \cite{BM_reb}.) Our proof is completely
different.  A direct calculation shows that $F_1=0$; thus, by
Theorem \ref{E}, $\gamma_1(\nu)\ge\lambda_1^D=-2\pi^2$, for all
$\nu$. In a preprint version of this paper, we made the conjecture that $\gamma_1(\nu)=\lambda_1^D=-2\pi^2$
for all $\nu$.
 This conjecture has   now been established by combining part (ii) with a very recent result    \cite{LLR07}
 which states  that in the case of $\frac12\frac{d^2}{dx^2}$ on an interval,
all of the eigenvalues are real, for all jump measures $\nu$.
\end{remark*}

The next result shows that
for the Laplacian on a square in $R^2$, one can
find a deterministic jump measure for
 which $\gamma_1(\nu)\neq\lambda_1^D$.

\begin{prop}\label{2dimBM}
 Consider the operator $\frac12\Delta$  in the square $(0,1)^2$. Then there exists a jump
 measure $\nu$ of the form $\delta_{x_0}$, for some $x_0\in (0,1)^2$,  for which
 $\gamma(\nu)>\lambda_1^D=-\frac{5\pi^2}2$.
\end{prop}



\begin{remark*}
Consider $L=\frac12\Delta$ in the cube $(0,1)^d$, $d\ge1$. The proof of Proposition
\ref{cube} shows that $F_1=0$; thus, by Theorem \ref{E}, $\gamma_1(\nu)\ge\lambda_1^D$.
Combining Proposition \ref{2dimBM} and Proposition \ref{cube}, it follows that for
$d\ge11$
or $d=2$,
there exists  a
jump measure $\nu$ for which $\gamma_1(\nu)>\lambda_1^D$.
Presumably, this holds for all $d\ge2$. Conversely, by Theorem \ref{thm:quasi inv}
 it follows that for all $d\ge1$ it is also always possible
to find a $\nu$ for which $\gamma_1(\nu)=\lambda_1^D$.
\end{remark*}

In all of the examples given so far, $\gamma_1(\nu)\ge\lambda_1^D$. The following result
shows that such is not always the case.
\begin{prop}\label{+-}
Let  $k_0=max\{n:\lambda_n^D=\lambda_1^D\}$. If $F_j\neq0$, for some $j\in\{1,2,\cdots,k_0\}$, then
it is always possible to find a jump measure $\nu$ for which $\gamma_1(\nu)>\lambda_1^D$
and it is always possible to find a jump measure $\nu$ for which
$\gamma_1(\nu)<\lambda_1^D$.
\end{prop}

\begin{remark*}
One can check that $L=\frac12\frac{d^2}{dx^2}+b\frac d{dx}$ on $(0,1)$, where
$b$ is a nonzero constant, is an example
where Proposition \ref{+-} is applicable.
\end{remark*}

We conclude this section with several open questions.

\bf\noindent Question 1.\rm\ In a preprint version of this paper, we asked whether all the eigenvalues
of $\mathcal L$ are real in
 the case that $L$ is
self-adjoint. It the very recent paper \cite{LLR07} it was shown that for $L=\frac12\Delta$ in a ball in $R^3$,
there exist deterministic jump measures, that is measures  of the form $\nu=\delta_x$, for which some
of the eigenvalues are not real. However, these non-real eigenvalues do not have maximal real part.
We still ask whether the eigenvalue with largest real part is real in the case
that $L$ is self-adjoint.

\noindent\bf Question 2.\rm\ Does the inequality $\gamma_1(\nu)<\lambda_0^D$ hold for all $\nu$
 when $L$ is self-adjoint? What about for general $L$?

\noindent \bf Question 3.\rm\ Does a lower bound exist for $\gamma_1(\nu)$ in terms of the
eigenvalues $\{\lambda_n^D\}_{n=0}^\infty$ of $L$?

\noindent \bf Question 4.\rm\ What can be said about the
continuity properties of $\gamma_1(\nu)$ as $\nu$ varies over
$\mathcal P(D)$, the space of probability measures on $D$ with the
topology of weak convergence?

\noindent \bf Remark.\rm\
 Note that if Question 1 is answered affirmatively, then Theorem \ref{upperbound} shows that the answer to Question
2 is affirmative in the case that $L$ is self-adjoint.

The proofs of the results stated in this section are grouped thematically and proved in the sections that follow.

\setcounter{equation}{0}\section{Proof of Theorem \ref{thm:quasi inv}}\rm\
 \rm
 We first prove the statement concerning the eigenvalues and eigenfunctions.
  Let $\phi^D_n$, $n\ge1$, denote an eigenfunction  for $L$ with the Dirichlet boundary condition,
 corresponding to the eigenvalue $\lambda_n^D$. Integration by parts yields
 $$
 \lambda_0^D \int_D \phi_n^D \witi{\phi^D_0} dx=
 \int_D \phi^D_n \witi{L}\witi{\phi^D_0} dx = \int_D L \phi^D_n \witi{\phi^D_0}dx =
  \lambda_n^D \int_D \phi_n^D \witi{\phi^D_0} dx,
  $$
 from which it follows  that $\int_D \phi^D_n \witi{\phi_0^D} dx =0$. Since $\phi_n^D|_{\partial D}=0$,
 it follows  that $\phi_n^D$ is in the domain of ${\mathcal L}$, and we conclude that
 $\lambda_n^D$ is an  eigenvalue for ${\mathcal L}$.
 Of course the function 1 is an eigenfunction for ${\mathcal L}$  corresponding to the
 eigenvalue 0.
Thus, to conclude the proof, it is enough to
show that the  equation
$L\psi=\lambda\psi$ with $\int_D\psi\witi{\phi^D_0}dx=\psi|_{\partial D}=c\neq 0$,
has a solution only when $\lambda=0$. Let  $\psi$ be a solution to the above equation.
Let $d\sigma$ denote the Lebesgue surface measure on $\partial D$.
Integrating by parts and using the boundary condition, we have
$$
\aligned&\lambda^D_0c=\lambda^D_0\int_D\witi{\phi^D_0}\psi dx=\int_D \psi \witi{L}\witi{\phi^D_0}dx\\
&=\int_D\witi{\phi^D_0}L\psi dx+\int_{\partial D}\frac12\psi a\nabla\witi{\phi^D_0}\cdot n d\sigma-\int_{\partial D}
\psi\witi{\phi^D_0}b\cdot nd\sigma\\
&=\lambda\int_D\witi{\phi^D_0}\psi dx+c\int_D\witi{L}\witi{\phi^D_0} dx=c\lambda+c\lambda^D_0,
\endaligned
$$
 where $n$ is the unit outward normal of $D$ at $\partial D$. Therefore  $\lambda=0$.

We now turn to the statement concerning the invariant measure.
 We denote by $G^D$ the operator from $L^1$ to $L^1$ given by $(G^D f) (x)= G^D(x,f)$. We denote its
adjoint on $L^{\infty}$ by $\witi{G^D}$. We have  $(\witi{G^D}g)(y)=G^D(g,y)$.
We now prove that $\witi{G^D}$ is compact.  For $\epsilon >0$, let $K_{\epsilon}$ denote the operator on $L^1$ defined by
 $$(K_{\epsilon} f)(y) = \int_{\epsilon}^{\epsilon^{-1}} p^D(s,f,y) ds.$$
 An argument similar but simpler than the one given in the proof of Lemma 1 in
 \cite{BP}, based on the  continuity of $p^D$ on $[\epsilon,\epsilon^{-1}]\times \ol{D}\times \ol{D}$, shows that $K_{\epsilon}$ is compact. Now,
 \begin{align}
 \nonumber
 \|\witi{G^D}f - K_{\epsilon}f\|_1 &= \int_D |\int_0^{\epsilon} p^D(s,f,y) ds  + \int_{\epsilon^{-1}}^{\infty} p^D(s,f,y) ds | dy \\
  \label{eq:witiG_approx}
 &\leq \epsilon \|f\|_1  + \|f\|_1 \int_{\epsilon^{-1}}^{\infty} \sup_{x\in  D} P_x^D(\tau_1 > s) ds.
\end{align}
 By \cite[Theorem 3.6.1]{pos_harm},
 $\lim_{s\to\infty} \frac 1s \log \sup_{x\in D}P_x(\tau_1 >s)=\lambda_0^D$. Thus, there exists a
 constant $C>0$ such that for all $\epsilon$ sufficiently small,
  $$\int_{\epsilon^{-1}}^{\infty} \sup_{x\in  D} P_x^D(\tau_1 > s) ds\leq C e^{\frac {\lambda_0^D}{2}\epsilon^{-1}}.$$
  Therefore, it follows from \eqref{eq:witiG_approx} that
 $$ \|\witi{G^D}-K_{\epsilon}\|_1 \leq \epsilon + C e^{\frac{\lambda_0^D}{2}\epsilon^{-1}}\underset{\epsilon \to 0}{\to }0.$$
 Consequently, $\witi{G^D}$ is compact.\par
 Assume now that $m\in {\mathcal P}$ is a fixed point for Inv. Since
 $d m(y) = \frac{G^D(m,y)}{G^D(m,\ch)}dy=\frac{\witi G^Dm}{G^D(m,\ch)}$, it follows that $m$ has density in $L^1$.
Therefore we may consider $m$ as an eigenfunction for $\witi{G^D}$,
corresponding  to the eigenvalue $\lambda \equiv G^D(m,\ch)$. Let $\phi^D_0>0$
denote the principal eigenfunction for $L$ with the Dirichlet
boundary condition corresponding to the eigenvalue $\lambda_0^D$.
Since $G^D=(-L)^{-1}$, we have
 $G^D \phi_0 = -(\lambda_0^D)^{-1}\phi_0^D$. Therefore,
 $$ \lambda \int_D m \phi^D_0 dx= \int_D \witi{G^D}m \phi^D_0dx = \int_D m G^D \phi^D_0 dx=
 -(\lambda_0^D)^{-1} \int_D m \phi^D_0dx.$$
 Since $m$ and $\phi^D_0$ are nonnegative,
 it follows that $\lambda = -(\lambda_0^D)^{-1}$.
Since $\witi G^D$ is compact,
 the Krein-Rutman theorem guarantees that
 the eigenvalue $-(\lambda_0^D)^{-1}$ for $\witi G^D$ is simple.
 Since $\witi{\phi_0}$ and $m$ are both
 eigenfunctions for $\witi{G^D}$ corresponding to the
 eigenvalue $-(\lambda_0^D)^{-1}$, and since
  $\int_D m dx = \int_D \witi{\phi_0} dx = 1$, we conclude that $m=\witi{\phi_0}$. \par
\medskip

\hfill $\square$

 \setcounter{equation}{0}\section{Proofs of Theorem
\ref{E}, Theorem \ref{real}, Theorem \ref{Ereversible}, Proposition \ref{cube}
and Proposition \ref{+-}}\rm\

\noindent\bf Proof of Theorem \ref{E}.\rm\ A number $\lambda\in C-\{0\}$ will be an eigenvalue
if and only there exists a function   $v$ satisfying $Lv=\lambda v$ and $v|_{\partial D}=\int_Dvd\nu$.
Let $u=v-c$,
where $c=v|_{\partial D}$. Then $u$ satisfies $Lu=\lambda u+K$,
where $K=\lambda c$.
On $L^2(D,\mu_{\text{rev}})$, the
function $u$ can be represented in the form
\begin{equation}\label{u}
u=\sum_{n=0}^\infty C_n\phi_n^D,
\end{equation}
for unknown constants $\{C_n\}_{n=0}^\infty$, and the
constant function 1 can be represented
by
\begin{equation}\label{1}
1=\sum_{n=0}^\infty F_n\phi_n^D,
\end{equation}
where $\{F_n\}_{n=0}^\infty$ is as in \eqref{FG}.
Since $u$ is a smooth function vanishing on $\partial D$, it is in the domain of
the self-adjoint operator $L$ acting on $L^2(D,\mu_{rev})$.
Thus,
from \eqref{u}, it follows that
\begin{equation}\label{Lu}
Lu=\sum_{n=0}^\infty\lambda_n^DC_n\phi_n^D.
\end{equation}
From \eqref{u}-\eqref{Lu} along with the fact that $Lu=\lambda u+K$, we obtain
\begin{equation}\label{key}
C_n\lambda_n^D=\lambda C_n+KF_n, \ n\ge0.
\end{equation}

We first show that the condition $E_\nu(\lambda)=0$ is necessary and sufficient for $\lambda\not\in\{\Lambda_n^D\}_{n=0}^\infty$ to
be an eigenvalue. Since we are now assuming that $\lambda$ is not in the spectrum of $L$, we may assume that $K\neq0$. Indeed,
if $K$ were equal to 0, then $v$ would vanish on $\partial D$ and consequently it would be an eigenfunction
for $L$. This would mean that $\lambda=\Lambda_n^D$, for some $n$.
From \eqref{key} we obtain
\begin{equation*}\label{Cnagain}
C_n=\frac{KF_n}{\lambda_n^D-\lambda},
\end{equation*}
and conclude that
\begin{equation}\label{uagain}
u=
\sum_{n=0}^\infty\frac{KF_n}{\lambda_n^D-\lambda}\phi_n^D.
\end{equation}

In order that $v$ be an eigenfunction, $v$ must satisfy
 $v|_{\partial D}=\int_Dv\ d\nu=c$. Since $u=v-c$, we require that $\int_Du\ d\nu=0$.
If $\nu$ has an $L^2(D,\mu_{rev})$-density, then using \eqref{uagain} and taking inner products shows that
$\int_Du\ d\nu=0$ if and only if
\begin{equation}\label{eigenvaluecondition}
\sum_{n=0}^\infty\frac{F_nG_n(\nu)}{\lambda_n^D-\lambda}=0,
\end{equation}
where $G_n(\nu)$ is as in \eqref{FG}.
Alternatively, if Assumption 1 holds,
then the formula for $u$ in \eqref{uagain} holds not only in $L^2(D,\mu_{rev})$, but also pointwise,
and from the bounded convergence theorem it follows again that $\int_Du \ d\nu=0$ if
and only if \eqref{eigenvaluecondition} holds.
We have thus shown  that the condition
$E_\nu(\lambda)=0$ is necessary and sufficient for a nonzero $\lambda\not\in\{\Lambda_n^D\}_{n=1}^\infty$
to be an eigenvalue. Furthermore, as the method uniquely specifies the corresponding eigenfunction
(up to a multiplicative constant), it follows that the multiplicity
of such an eigenvalue is 1.

We now consider the possibility that $\lambda=\Lambda^D_{n_0}$ is an eigenvalue, where $n_0$ is a nonnegative integer.
Let
$S_{n_0}$ denote the $d_{n_0}$-dimensional eigenspace corresponding to the eigenvalue
$\Lambda^D_{n_0}$ of $L$.
Let $S^G_{n_0}(\nu)=\{w\in S_{n_0}: \int_D wd\nu=0\}$ and let $S^F_{n_0}=\{w\in S_n:\int_Dwd\mu_{\text{rev}}=0\}$.
Clearly, each of these latter two spaces is either $(d_{n_0}-1)$-dimensional or $d_{n_0}$-dimensional.

Consider first the case that $S_{n_0}^F$ is $(d_{n_0}-1)$-dimensional.
There exists an $m_0$ such that $\lambda_{m_0}^D=\Lambda_{n_0}^D$ and $F_{m_0}=\int_D\phi_{m_0}d\mu_{\text{rev}}\neq0$.
But then \eqref{key} will hold with $n=m_0$ and $\lambda=\Lambda_{n_0}$ if and only if $K=0$.
But if $K=0$, then $u=v$, $v|_{\partial D}=\int_Dvd\nu=0$ and $\Delta v=\Lambda^D_{n_0}v$.
Thus,   $v$ belongs
to $S^G_{n_0}(\mu)$. Consequently, the multiplicity of $\Lambda^D_{n_0}$ will be either
$d_{n_0}-1$ or $d_{n_0}$, depending on which of these numbers is the dimension of $S^G_{n_0}(\mu)$.
In particular, if $n_0=0$, then $d_{n_0}=1$ and $S_{n_0}^F=S_{n_0}^G(\nu)=\{0\}$ since $\phi_0^D>0$. Thus,
$\lambda_0^D=\Lambda_0^D$ can never be an eigenvalue.

Now consider the case that $S^F_{n_0}$ is $d_{n_0}$-dimensional.
In this case, $F_m=0$, for all $m$ such that $\lambda_m^D=\Lambda^D_{n_0}$.
We first look for eigenfunctions for which $K\neq0$. Solving \eqref{key}
gives
\begin{equation*}
\begin{cases}
C_n=\frac{KF_n}{\lambda^D_n-\lambda},\ \text{for all}\ n\ \text{such that}\ \lambda^D_n\neq\Lambda_{n_0}^D;\\
C_n\ \text{is arbitrary}, \text{for all}\ n\ \text{such that}\ \lambda^D_n=\Lambda_{n_0}^D.
\end{cases}
\end{equation*}
Writing $C_n=Kc_n$, for $n$ such that $\lambda^D_n=\Lambda_{n_0}^D$,
and employing the
 same reasoning as in \eqref{uagain} and \eqref{eigenvaluecondition} yields
\begin{equation}\label{kn}
\sum_{n:\lambda_n^D\neq\Lambda_{n_0}^D} \frac{ F_nG_n(\nu)}{\lambda_n^D-\Lambda_{n_0}^D}
+\sum_{n:\lambda_n^D=\Lambda_{n_0}^D}c_nG_n(\nu)=0.
\end{equation}
There are two cases to consider---when $S^G_{n_0}(\nu)$ is $(d_{n_0}-1)$-dimensional and when it is $d_{n_0}$-dimensional.
In the latter case, $G_n(\nu)=0$, for all $n$ satisfying $\lambda_n^D=\Lambda^D_{n_0}$.
Thus,  \eqref{kn} reduces to $E_\nu(\Lambda_{n_0}^D)=0$.
If this equation is satisfied, we obtain one eigenfunction with $K\neq0$, and if it is not satisfied, we obtain
no such eigenfunctions. Since
$S^G_{n_0}(\nu)$ is $d_{n_0}$-dimensional, there are also  $d_{n_0}$
additional linearly independent eigenfunctions with $K=0$.
Thus,  the multiplicity is either $d_{n_0}+1$ or $d_{n_0}$, depending on whether
or not $E_{\nu}(\Lambda_{n_0}^D)=0$.

Now consider the case that
$S^G_{n_0}(\nu)$ is $(d_{n_0}-1)$-dimensional.
Since we may  choose the orthonormal basis $\{\phi^D_m\}_{\{m:\lambda_m^D=\Lambda_{n_0}^D\}}$
corresponding to the eigenspace $S_{n_0}$ however we like, we may assume
without loss of generality, that
$G_m(\nu)=\int_D\phi_md\nu=0$,  for all but one of the $m$ for which $\lambda_m^D=\Lambda_{n_0}^D$.
Denote the single $m$ for which this is not true by $m_0$.
Then \eqref{kn} reduces to
\begin{equation*}
\sum_{n:\lambda_n^D\neq\Lambda_{n_0}^D} \frac{ F_nG_n(\nu)}{\lambda_n^D-\Lambda_{n_0}^D}
+c_{m_0}G_{m_0}(\nu)=0.
\end{equation*}
The above equation is uniquely solvable
for $c_{m_0}$, and thus yields one eigenfunction with $K\neq0$.
Since $S^G_{n_0}(\nu)$ is $(d_{n_0}-1)$-dimensional, there are also  $d_{n_0}-1$
additional linearly independent eigenfunctions with $K=0$; thus the multiplicity
is $d_{n_0}$.
\hfill $\square$
\medskip

\noindent \bf Proof of Theorem \ref{real}.\rm\ By Theorem \ref{E},
a complex number $\lambda=\alpha+i\beta$, with $\beta\neq0$ will
be an eigenvalue for $\mathcal L$ if and only if
\[
\sum_{n=0}^\infty\frac{F_nG_n(\nu)}{\lambda_n^D-\lambda}=0.
\]
We can rewrite this as
\begin{equation}\label{complex1}
\begin{aligned}
&\sum_{n=0}^\infty\frac{F_nG_n(\nu)(\lambda_n^D-\alpha)}{(\lambda_n^D-\alpha)^2+\beta^2}=0;\\
&\sum_{n=0}^\infty\frac{F_nG_n(\nu)\beta}{(\lambda_n^D-\alpha)^2+\beta^2}=0.
\end{aligned}
\end{equation}
Clearly, the two equations in \eqref{complex1}
hold if and only if the following two equations hold:
\begin{equation}\label{complex2}
\begin{aligned}
&\sum_{n=0}^\infty\frac{F_nG_n(\nu)\lambda_n^D}{(\lambda_n^D-\alpha)^2+\beta^2}=0;\\
&\sum_{n=0}^\infty\frac{F_nG_n(\nu)}{(\lambda_n^D-\alpha)^2+\beta^2}=0.
\end{aligned}
\end{equation}
Since $F_0G_0(\nu)$ is always positive,
neither equation in \eqref{complex2} can hold if $F_nG_n(\nu)\ge0$, for all $n\ge1$.
Consider now either the case that  $F_nG_n(\nu)\le0$, for all $n\ge1$,
or alternatively, the case that  $F_nG_n(\nu)$ is nonzero for no more than
two values of $n$.
Rewriting \eqref{complex2} as
$$
\begin{aligned}
&
\frac{F_0G_0(\nu)}{(\lambda_0^D-\alpha)^2+\beta^2}+
\sum_{n=1}^\infty\frac{F_nG_n(\nu)}{(\lambda_n^D-\alpha)^2+\beta^2}\frac{\lambda_n^D}{\lambda_0^D}=0;\\
&\frac{F_0G_0(\nu)}{(\lambda_0^D-\alpha)^2+\beta^2}+
\sum_{n=1}^\infty\frac{F_nG_n(\nu)}{(\lambda_n^D-\alpha)^2+\beta^2}=0,
\end{aligned}
$$
it follows that the two equations in \eqref{complex2}
cannot hold simultaneously.
\hfill $\square$

\noindent \bf Proof of Theorem \ref{Ereversible}.\rm\ \it i.\rm\
Since $G_n(\mu_{rev})=F_n$, it follows from part (i) of Theorem
\ref{real} that all the eigenvalues of $\mathcal L$ are real.

\it \noindent ii.\rm\ By part (i), $\gamma_1=\gamma_1(\mu_{\text{rev}})$ is itself  an eigenvalue; let $\phi_1$ denote
a corresponding eigenfunction. Let $\psi_1=\phi_1-c$,
where $c=\phi_1|_{\partial D}=\int_D\phi_1d\mu_{\text{rev}}$.
Then  $\psi_1|_{\partial D}=\int_D\psi_1d\mu_{\text{rev}}=0$ and
 $L\psi_1=\gamma_1\psi_1+\gamma_1c$.
 Multiplying this equation
by $\psi_1\exp(2Q)$ and integrating by parts gives
$$
\gamma_1=-\frac{\frac12\int_D(\nabla \psi_1a\nabla \psi_1)d\mu_{\text{rev}}}{\int_D\psi_1^2d\mu_{\text{rev}}}.
$$

On the other hand, consider the quotient
$\frac{\frac12\int_D(\nabla ua\nabla u)d\mu_{\text{rev}}}{\int_Du^2d\mu_{\text{rev}}}$.
By standard methods, the infimum
of this quotient over functions  $0\neq u\in H_0^1(D)$ satisfying $u|_{\partial D}=\int_Dud\mu_{\text{rev}}=0$ exists.
We denote this infimum by $-\Gamma>0$.
To identify the minimum, we use a Lagrange multiplier and vary the quantity
$\frac12\int_D(\nabla ua\nabla u)d\mu_{\text{rev}}+k\int_Du^2d\mu_{\text{rev}}$ over
functions $u$ satisfying the above restriction, where $k$ is a free parameter.
 A minimizer $\psi$ must satisfy the equation
 $\int_Dq(L\psi-k\psi)  ~d\mu_{\text{rev}}=0$, for all $q$ satisfying the above restriction.
 From this one concludes that $L\psi=k\psi+C$, for some constant $C$. Multiplying this
equation by $\psi$, integrating both sides with respect to $d\mu_{\text{rev}}$,
 and integrating by parts, one finds that $k=\Gamma$.
Letting $\phi=\psi+\frac C\Gamma$, it follows  that $\phi$ satisfies
$L\phi=\Gamma\phi$ and $\phi|_{\partial D}=\int_D\phi d\mu_{\text{rev}}$.

\noindent\it iii.\rm\
By part (i) and  the definition of
  $\gamma_1(\mu_{\text{rev}})$, it follows that
$\gamma_1(\mu_{\text{rev}})$
   is the largest nonzero eigenvalue of $\mathcal L$.
And then by Theorem \ref{upperbound} it follows that $\gamma_1(\mu_{rev})<\lambda_0^D$.
In Theorem \ref{E}, note  that when $\nu=\mu_{rev}$, then $G_n(\mu_{rev})=F_n$. Consequently
$E_{\mu_{rev}}(\lambda)=\sum_{n=0}^\infty\frac{F_n^2}{\lambda_n^D-\lambda}$.
Since $E_{\mu_{rev}}(\lambda)$ is continuous for $\lambda\in(\lambda_1^D,\lambda_0^D)$,
 since $E_{\mu_{rev}}((\lambda_0^D)^-)=\infty$
and since $E_{\mu_{rev}}((\lambda_1^D)^+)=-\infty$ holds if $F_j=\int_D\phi_j^D d\mu_{\text{rev}}\neq0$,
for some $j\in\{1,2,\cdots, k_0\}$,
it follows that $E_{\mu_{rev}}$ possesses a root in $(\lambda_1^D,\lambda_0^D)$
 if $F_j\neq0$ for some $j\in\{1,2,\cdots,k_0\}$.
It now follows from
  Theorem \ref{E}
 that $\lambda_1^D\le\gamma_1(\mu_{rev})<\lambda_0^D$,
 with strict inequality if $F_j\neq0$, for some $j\in\{1,2,\cdots,k_0\}$.
Furthermore, since $E_{\mu_{rev}}$ is increasing on $(\lambda_1^D,\lambda_0^D)$,
if follows that in the case that $F_j=0$ for all $j\in\{1,2,\cdots,k_0\}$, the strict
inequality will hold if and only if $E_{\mu_{rev}}(\lambda_1^D)<0$. This inequality
can be rewritten as \eqref{borderline}.
\hfill $\square$

\medskip

\noindent \bf Proof of Proposition \ref{cube}.\rm\
By Theorem \ref{Ereversible}-iii, $\gamma_1(\mu_{rev})<\lambda_0^D$. To prove the rest of the proposition,
we apply \eqref{borderline} from Theorem \ref{Ereversible}. The complete, orthonormal sequence of
eigenfunctions on $L^2(D,l_d)$ for $\frac12\Delta$ on $D\equiv(0,1)^d$
with the Dirichlet boundary condition is given by \linebreak $\{2^{\frac
d2}\prod_{j=1}^d\sin n_j\pi x_j\}_{n_1,\cdots,n_d=1}^\infty$. The
corresponding eigenvalues are\linebreak
$\{-\frac{\pi^2}2\sum_{j=1}^dn_j^2\}_{n_1,\cdots,n_d=1}^\infty$. We
will denote these eigenfunctions and eigenvalues respectively by
$\phi^D_{n_1,\cdots, n_d}$ and $\lambda^D_{n_1,\cdots, n_d}$. We
have
\begin{equation}\label{F}
F_{n_1,\cdots,n_d}\equiv\int_D\phi^D_{n_1,\cdots,n_d}\ dx=
\begin{cases}
\frac{2^{\frac32d}}{\pi^d\prod_{j=1}^dn_j}, \ \text{if}\ n_j\
\text{is odd for all}\ j;\\ 0, \ \text{otherwise}.
\end{cases}
\end{equation}
In the present context, the terms $F_0,\lambda_0^D$ and $\lambda_1^D$
appearing in Theorem
\ref{Ereversible} are given respectively by
$F_{1,\cdots,1}=\frac{2^{\frac32d}}{\pi^d}$,
$\lambda_{1,\cdots,1}=-\frac{d\pi^2}2$ and $\lambda_{n_1,\cdots,n_d}=-\frac{(d+3)\pi^2}2$,
where $(n_1\cdots,n_d)$ satisfies $\sum_{j=1}^d n_j=d+1$. From \eqref{F},
we have $F_{n_1,\cdots,n_d}=\int_D\phi_{n_1,\cdots,n_d}^D~d\mu_{\text{rev}}=0$, if
$\sum_{j=1}^dn_j=d+1$. Thus, \eqref{borderline} is applicable.

Using $\{\lambda_{n_1,\cdots,n_d}\}$ and $\{F_{n_1,\cdots,n_d}\}$ in place
of the labeling $\{\lambda_n\}$ and $\{F_n\}$ in
the inequality \eqref{borderline}, we find  that after cancellations
the inequality can be written
as
\begin{equation}\label{condition}
\sum_{\substack {n_1,\cdots,n_d\ \text{odd}\\ (n_1,\cdots,n_d)\neq (1,\cdots,1)}}
\frac1{\prod_{j=1}^dn^2_j\left(\sum_{j=1}^dn_j^2-d-3\right)}>\frac13.
\end{equation}
Thus, by \eqref{borderline}, \eqref{condition} is a necessary and sufficient
condition in order that $\gamma_1(\mu_{\text{rev}})>\lambda_1^D$,
and if the condition does not hold, then $\gamma_1(\mu_{\text{rev}})=\lambda_1^D$.
Denote the left hand side  of \eqref{condition} by $H_d$.
If one considers $H_{d+1}$,  but restricts the summation to those
multi-indices $(n_1,\cdots,n_{d+1})$ for which $n_{d+1}=1$, the resulting
quantity is $H_d$. Thus the left hand side of \eqref{condition} is monotone increasing in
$d$.
A direct calculation shows that the inequality in \eqref{condition} does not hold if $d=1$.
On the other hand, by considering
 the contribution to the left hand side of \eqref{condition} only
from those multi-indices satisfying $\sum_{j=1}^dn_j=d+2$, it is easy to check that
the inequality in \eqref{condition} holds if $d\ge15$.
From these observations we conclude that there exists a $d^*\in[2,15]$ such that
$\gamma_1(\mu_{\text{rev}})>\lambda_1^D$, if
$d\ge d^*$, and $\gamma_1(\mu_{\text{rev}})=\lambda_1^D$, if $d<d^*$.
Numerical calculations shows
that in fact $d^*=11$.
\hfill $\square$



\medskip
\noindent \bf Proof of Proposition \ref{+-}.\rm\
Without loss of generality, assume that $F_1>0$.
Choose $\nu^\pm$ with density $\nu^\pm(x)=c_\pm(\phi_0^D(x)\pm\epsilon \phi_1^D(x))$,
where $\epsilon>0$ is sufficiently small so that $\nu^\pm(x)\ge0$, for all $x\in D$, and
where $c_\pm>0$ is a normalizing constant so that $\nu^\pm$ is a probability density.
(It is possible to choose such an $\epsilon>0$ because the Hopf
maximum principle guarantees that the normal derivative
$\nabla\phi_0^D(x)\cdot n\neq0$, for $x\in\partial D$.)
Recall that $G_n(\nu^\pm)=\int_D\phi_n^D \ d\nu^\pm$.
Thus, $G_0(\nu^\pm)=c_\pm$, $G_1(\nu^\pm)=\pm c_\pm\epsilon$
and $G_n(\nu^\pm)=0$, for $n\ge2$.
From the definition of $E_{\nu^\pm}$ in Theorem \ref{E}, we have
$E_{\nu^\pm}(\lambda)=\frac{c_\pm F_0}{\lambda_0^D-\lambda}\pm\frac{c_\pm\epsilon F_1}{\lambda_1^D-\lambda}$.
Thus, $E_{\nu^-}(\lambda)\neq0$ for $\lambda\in(\lambda_1^D,\lambda_0^D)$. By Theorems \ref{real}
and \ref{upperbound}, $\gamma_1(\nu^-)<\lambda_0^D$. Thus, we conclude from
Theorem \ref{E} that $\gamma_1(\nu^-)<\lambda_1^D$.
Since $E_{\nu^+}((\lambda_0^D)^-)=\infty$
and $E_{\nu^+}((\lambda_1^D)^+)=-\infty$,
$E_{\nu^+}(\lambda)$ possesses a root in $(\lambda_1^D,\lambda_0^D)$. Thus, by Theorem \ref{E}.
$\gamma_1(\nu^+)>\lambda_1^D$.
\hfill $\square$

 \setcounter{equation}{0}\section{Proof of Theorem \ref{upperbound}}\rm\
By assumption, $\gamma_1(\nu)$ is a real eigenvalue for $\mathcal
L$. We need to show that $\gamma_1(\nu)<\lambda_0^D$. Let $u$
denote a corresponding eigenfunction, and let $c=u|_D=\int_Du
d\nu$. We first show that $\gamma_1(\nu)\neq\lambda_0^D$. Assume
to the contrary. In this case, $c\neq0$.  Indeed, otherwise
$\phi_0^D$ and $u$ would both be eigenfunctions for the principal
eigenvalue $\lambda_0^D$ of $L$ with the Dirichlet boundary
condition. Furthermore, $\phi_0^D$ and $u$ would be linearly
independent since $\phi_0^D$ does not change sign, whereas
$\int_Dud\nu=c=0$. This would then contradict the simplicity of
the principal eigenvalue $\lambda_0^D$. Integrating by parts
twice, exploiting the form of the reversible operator $L$ and the
reversible measure, we have
\[
\begin{aligned}
&\int_D\phi_0^DLu\ d\mu_{rev}=\int_{\partial D}\phi_0^Da\nabla u\cdot n
\exp(2Q)d\sigma\\
&-\int_{\partial D}ua\nabla \phi_0^D\cdot n\exp(2Q)d\sigma+\int_DuL\phi_0^D\ d\mu_{rev},
\end{aligned}
\]
which reduces to
\begin{equation}\label{0}
\int_{\partial D}ua\nabla \phi_0^D\cdot n\exp(2Q)d\sigma=0.
\end{equation}
However,
\begin{equation}\label{00}
\begin{aligned}
&\int_{\partial D}ua\nabla \phi_0^D\cdot n\exp(2Q)d\sigma=c\int_{\partial D}
a\nabla \phi_0^D\cdot n\exp(2Q)d\sigma\\
&=c\int_DL\phi_0^Dd\mu_{rev}=c\lambda_0^D\int_D\phi_0^Dd\mu_{rev}.
\end{aligned}
\end{equation}
Now \eqref{0} and \eqref{00} give $c=0$, which is a contradiction.

Now we show that $\gamma_1(\nu)\ngtr\lambda_0^D$. Assume to the contrary.
By the Feynman-Kac formula,
$u(Y(t\wedge\tau_D))\exp(-\gamma_1(\nu)(t\wedge\tau_D))$ is a martingale.
Thus,
\begin{equation}\label{martingale}
\ex^D_xu(Y(t\wedge\tau_D))\exp(-\gamma_1(\nu)(t\wedge\tau_D))=u(x).
\end{equation}
Since $\gamma_1(\nu)>\lambda_0^D$, we have
$\ex^D_x\exp(-\gamma_1(\nu)\tau_D)<\infty$ \cite[chapter 3]{pos_harm}. Thus letting $t\to\infty$ in
\eqref{martingale} and applying the dominated convergence theorem gives
\begin{equation}\label{u(x)}
u(x)=c\ex^D_x\exp(-\gamma_1(\nu)\tau_D).
\end{equation}
It follows from \eqref{u(x)} that $c\neq0$.
Integrating both sides of \eqref{u(x)} against $\nu$ now gives
\begin{equation}\label{contradiction}
\ex^D_\nu\exp(-\gamma_1(\nu)\tau_D)=1,
\end{equation}
which is a contradiction.\hfill$\square$

\begin{remark*}
If one does not assume that the nonzero eigenvalue with largest real part is real,
 the calculation in the above proof
can be made with $\gamma_1(\nu)$ replaced by $\lambda_1(\nu)$,
where $\lambda_1(\nu)$ is an eigenvalue for $\mathcal L$ whose
real part is $\gamma_1(\nu)$. One arrives at \eqref{contradiction}
with $\gamma_1(\nu)$ replaced by $\lambda_1(\nu)$. However, since
$\lambda_1(\nu)$ can be complex-valued, \eqref{contradiction} no
longer constitutes a contradiction.
\end{remark*}

\setcounter{equation}{0}\section{Proof of Proposition \ref{1dimBM}}
\noindent \it i.\rm\
 The eigenvalue problem  for ${\mathcal L}$ is
 \begin{equation}
\label{BM8}
 \begin{cases}
  &\frac 12 u''  = \lambda u\
  \mbox{in }(0,1);\\ & u(p) = u(0)=u(1).\end{cases}
\end{equation}
Every solution $u$ to the differential equation in \eqref{BM8} is given by
 $$u(x)=A \cos \kappa x + B\sin \kappa x,$$
 where $\lambda = - \frac 12 \kappa^2$.
 In order that such a solution also
 satisfy the boundary condition in \eqref{BM8},
  the following system of linear equations must have a nontrivial solution:
\begin{align}
\nonumber
 A\left(1- \cos \kappa p \right) - &B\sin \kappa p  =0; \\
\label{eq:BM8_det}
 A\left(1  - \cos \kappa  \right)  - &B\sin  \kappa =0.
\end{align}
The determinant of the linear system above is
\begin{equation}\label{sys}
 \begin{aligned}
 -(1-\cos \kappa p)\sin \kappa &+ \sin \kappa p (1-\cos \kappa)\\
 & =
 \sin \kappa (1-p)-\sin \kappa +\sin \kappa p \\
 & =
 2\sin \frac{\kappa (1-p)}{2}\cos \frac{\kappa (1-p)}{2} +2 \cos \frac{\kappa (p+1)}{2}\sin {\kappa(p-1)}{2} \\
& =
  2\sin \frac{\kappa (1-p)}{2}\left( \cos \frac{\kappa (1-p)}{2} - \cos \frac{\kappa (p+1)}{2}\right)\\
 &=
 4\sin \frac{\kappa (1-p)}{2}\sin \frac{\kappa}{2}\sin \frac{\kappa p}{2}.
\end{aligned}
\end{equation}
 Since $\sin x = 0$ if and only if $x = \pi n$ for some integer $n$, the solutions $\kappa$ of \eqref{sys} are all real and are given by
 $$ \frac {2\pi n}{1-p},~{2\pi n}, \frac{2 \pi n}{p},~n\in \Z.$$
 Therefore, the eigenvalues for ${\mathcal L}$ are
 $$ - \frac {2 \pi^2 n^2}{(1-p)^2},~- 2 \pi^2 n^2,~-\frac {2\pi^2 n^2}{p^2},~n\in\N.$$
 Thus, the non-zero eigenvalue with maximal real part is $-2 \pi^2$.

\noindent \it ii.\rm\
By assumption, the nonzero eigenvalue with largest real part is real,
and we know that it is negative.
Thus, $\gamma_1(\nu)$ is the largest negative
number $\gamma$ for which there is a solution to the following problem:
 \begin{equation}
\label{eq:BM8_ev_prob}
 \begin{cases}
  \frac12u''= \gamma u&\mbox{in }(0,1);\\u(0)=u(1)= \int_0^1 u\ d\nu.\end{cases}
\end{equation}
 Every solution to the differential equation in \eqref{eq:BM8_ev_prob} with $\gamma<0$ is of
 the form
 $$
 u(x)=A \cos \kappa x + B\sin \kappa x,
 $$
 where $\gamma = - \frac 12 \kappa^2$, for some $\kappa\in R-\{0\}$.
 In order that \eqref{eq:BM8_ev_prob} have a solution, the following system of linear equations must have a nontrivial solution:
\begin{align}
\nonumber
 A\left(1- \int_0^1\cos \kappa x \ d\nu\right) - &B\int_0^1\sin \kappa x\ d\nu=0; \\
\label{eq:BM8_det}
 A\left(1  - \cos \kappa  \right)  - &B\sin  \kappa =0.
\end{align}
The determinant of the linear system above is
\[
 \begin{aligned}
& -(1-\int_0^1\cos \kappa x\ d\nu)\sin \kappa + (1-\cos \kappa) \int_0^1\sin \kappa x\ d\nu\\
& =
\int_0^1 \sin \kappa (1-x)\ d\nu-\sin \kappa +\int_0^1\sin \kappa x\ d\nu \\
&  =
 2\int_0^1\sin \frac{\kappa (1-x)}{2}\cos \frac{\kappa (1-x)}{2}\ d\nu +
 2\int_0^1 \cos \frac{\kappa (x+1)}{2}\sin \frac{\kappa(x-1)}{2}\ d\nu \\
& =
  2\int_0^1\sin \frac{\kappa (1-x)}{2}\left( \cos \frac{\kappa (1-x)}{2} - \cos \frac{\kappa (x+1)}{2}\right)d \nu\\
 &=
 4\int_0^1\sin \frac{\kappa (1-x)}{2}\sin \frac{\kappa}{2}\sin \frac{\kappa x}{2}\ d\nu.
\end{aligned}
\]
Note that
$\sin \frac{\kappa (1-x)}{2}\sin \frac{\kappa}{2}\sin \frac{\kappa x}{2}>0$,
for $\kappa\in(0,2\pi)$ and $x\in(0,1)$, while
the reverse inequality holds for $\kappa\in(-2\pi,0)$ and $x\in(0,1)$. Thus, it follows that
\eqref{eq:BM8_det} has no solution $\kappa\in(-2\pi,2\pi)-\{0\}$. On the other
hand,
$\sin \frac{\kappa (1-x)}{2}\sin \frac{\kappa}{2}\sin \frac{\kappa x}{2}\equiv0$,
for $\kappa=\pm2\pi$. Thus, it follows that $\gamma=-\frac12(2\pi)^2=-2\pi^2$ is the largest real
nonzero solution to \eqref{eq:BM8_ev_prob}.
\hfill $\square$.

 \setcounter{equation}{0}\section{Proof of Proposition \ref{2dimBM}}\rm\
By Theorem \ref{conditions} (part (ii) or part(iii)), Theorem \ref{E}  holds for all jump measures $\nu$.
We will show that for an appropriate $\nu$ one has $E_\nu(\lambda)=0$, for some $\lambda>\lambda_1^D$.

We use the notation and the calculations in the proof of
Proposition \ref{cube}. Let  $\nu_0=\delta_{(\frac19,\frac19)}$.
We have $G_{n_1,n_2}(\nu_0)=\int_D\phi_{n_1,n_2}^D\
d\nu_0=\phi_{n_1,n_2}^D(\frac19,\frac19)=
2\sin(\frac{n_1}9\pi)\sin(\frac{n_2}9\pi)$. Then from the
definition of $E_{\nu_0}(\lambda)$ it follows that
\begin{equation}\label{Enu0}
\begin{aligned}
&E_{\nu_0}(\lambda)=\\
&C\sum_{m_1,m_2=0}^\infty\frac{\sin(\frac{2m_1+1}9\pi)
\sin(\frac{2m_2+1}9\pi)}{(2m_1+1)(2m_2+1)\left((2m_1+1)^2+(2m_2+1)^2+\frac2{\pi^2}\lambda\right)},
\end{aligned}
\end{equation}
for an appropriate negative constant $C$. We will  show that the
equation $E_{\nu_0}(\lambda)=0$ has a root
$\lambda\in(\lambda_1^D,\lambda_0^D)=(-\frac52\pi^2,-\pi^2)$.
Note that $E_{\nu_0}((-\pi^2)^-)=\infty$.
Thus, it suffices to show that $E_{\nu_0}(-\frac{5\pi^2}2)<0$.
This can be checked using a program such as Mathematica, or alternatively, by a page and a half of estimates
which we refrain from reproducing here.

 \setcounter{equation}{0}\section{Proof of Theorem
\ref{conditions}}\rm\ For all three parts of the theorem, we will need
the following comparison result. By the mini-max principle
\cite{RS78}, one can compare the eigenvalues
$\{\lambda_n^D\}_{n=0}^\infty$ of $L$ in $D$ to those of
$\frac12\Delta$
in $(0,1)^d$, and conclude that there exist
$c_1,c_2>0$ (depending on $L$ and $D$) such that $c_1\hat \lambda_n\le \lambda_n^D\le c_2\hat \lambda_n$, where
$\{\hat\lambda_n\}_{n=0}^\infty$ are the eigenvalues for $\frac12\Delta$ in $(0,1)^d$, labelled in nonincreasing order.

\it\noindent i.\rm\ It is known that the eigenfunctions $\{\phi_n^D\}_{n=0}^\infty$ are uniformly
bounded \cite[pp.270-273]{Ince56}. Thus,
using the Cauchy-Schwarz inequality, it is enough to show that
$\sum_{n=0}^\infty\frac1{(\lambda_n^D)^2}<\infty$.
By the comparison principle above, it suffices to show the above inequality
in the case that $L=\frac12\frac{d^2}{dx^2}$ in $D=(0,1)$. In this case,
$\lambda_n^D=-\frac{(n+1)^2\pi^2}2$.

\it \noindent ii.\rm\ When $L$ is the Laplacian on a Riemannian manifold, it is known that
 $|\phi_n^D|\le C|\lambda_n^D|^{\frac{d-1}4}$, for some $C>0$ \cite{Gri02}.
 Using this and applying the Cauchy-Schwarz inequality,
 it follows that   $\sum_{n=0}^\infty \frac{F_n}{\lambda_n^D}\phi_n^D(x)$,
 converges uniformly and absolutely when $d=2$ if
 $\sum_{n=0}^\infty|\lambda_n^D|^{-\frac32}$ converges.
By the mini-max principle above, it suffices to show the above
inequality in the case that $L=\frac12\Delta$ on $(0,1)^2$.
But this then follows from Weyl's asymptotic distribution of eigenvalues \cite{RS78} which gives
$\lambda_n\sim cn$.

\it \noindent iii.\rm\ As in part (i),  it is enough to show that
$\sum_{n=0}^\infty\frac1{(\lambda_n^D)^2}<\infty$, and by the comparison
principle above, it suffices to show the above inequality in the case that
$L=\frac12\Delta$ on $(0,1)^d$, $d\le3$.
But this then follows from Weyl's asymptotic distribution of eigenvalues \cite{RS78}, which gives
$\lambda_n\sim cn^{\frac2d}$.

\hfill $\square$


\begin{thebibliography}{99}




\bibitem{BP} Ben Ari, I. and Pinsky, R. G., Ergodic behavior of diffusions with random jumps from the boundary,
submitted.


\bibitem{Gri02} Grieser, D., Uniform bounds for eigenfunctions of the Laplacian on manifolds with boundary,
\emph{Comm. Partial Differential Equations} \textbf{27}
      (2002), 1283-1299.

\bibitem{BM_8} Grigorescu, I.  and Kang, M.,
     Brownian motion on the figure eight,
   \emph{J. Theoret. Probab.} \textbf{15}, (2002), 817--844.
      


\bibitem{BM_reb}
  Grigorescu, I. and  Kang, M.,
  Ergodic properties of multidimensional Brownian motion with rebirth,
  preprint, url:\ http://www.math.miami.edu/~igrigore/pp/gn.pdf.

\bibitem{CJ} Carslaw, J. and Jaeger, J. Conduction of Heat in solids,
\emph{Reprint of the second edition. Oxford Science Publications. The Clarendon Press, Oxford University Press, New York}, (1988).


\bibitem{Ince56} Ince, E. L.,
   Ordinary Differential Equations,
   \emph{Dover},
   New York,
     (1956).


\bibitem{LLR07}
Leung, Y., Li, W. and Rakesh,
Spectral analysis of Brownian motion with jump boundary,
\emph{preprint}


\bibitem{Pin85} Pinsky, R. G.,
     On the convergence of diffusion processes conditioned to remain
     in a bounded region for large time to limiting positive recurrent diffusion processes,
   \emph{Annals of Probab.} \textbf{13} (1985), 363-378.

\bibitem{pos_harm} Pinsky, R. G.,
     Positive Harmonic Functions and Diffusion,
    \emph{Cambridge Studies in Advanced Mathematics}
    \textbf{45}, \emph{Cambridge University Press},
      (1995)

\bibitem{RS78}
    Reed, M.  and Simon, B.,
     Methods of Modern Mathematical Physics, IV,  Analysis of
              Operators,
 \emph{Academic Press [Harcourt Brace Jovanovich Publishers]},
   New York, (1978).

\bibitem{vandenberg} van den Berg, M. Heat content of a Riemannian manifold with a perfect conducting boundary,
\emph{Potential Analysis} \textbf{19}, (2003), 89-98.



\end{thebibliography}
\end{document}